\documentclass[11pt, reqno, english]{smfart}

\usepackage[T1]{fontenc}
\usepackage[english,francais]{babel}

\usepackage{amssymb, url, xspace, smfthm}

\usepackage{graphicx}
\usepackage{subfigure}
\usepackage{psfrag}

\usepackage[svgnames, table]{xcolor}
\definecolor{oneblue}{rgb}{0.0, 0.0, 0.85}

\usepackage[usenames, dvipsnames]{pstricks}
\usepackage{epsfig}   
\usepackage{pst-grad} 
\usepackage{pst-plot} 

\usepackage[colorlinks,
            urlcolor=oneblue,
            linkcolor=LightSlateGrey,
            citecolor=gray,
            bookmarksopen=false,
            pagebackref]{hyperref}
\hypersetup{pdfpagemode=UseNone}

\addtocontents{toc}{\setcounter{tocdepth}{2}}

\makeatletter
    \def\ps@copyright{\ps@empty
    \def\@oddfoot{\hfil\small\copyright NumHyp--2015}}
\makeatother

\usepackage{acronym}
\acrodef{SGN}[SGN]{Serre--Green--Naghdi}
\acrodef{eSGN}[eSGN]{extended Serre--Green--Naghdi}

\newcommand{\BibTeX}{{\scshape Bib}\kern-.08em\TeX}
\newcommand{\T}{\S\kern .15em\relax }
\newcommand{\AMS}{$\mathcal{A}$\kern-.1667em\lower.5ex\hbox
        {$\mathcal{M}$}\kern-.125em$\mathcal{S}$}

\def\R{\mathbb{R}}  
\def\Dt{\partial_t}
\def\Dx{\partial_x}
\def\Dv{\partial_v}
\def\DX{{\Delta x}}
\def\DT{{\Delta t}}
\def\Nbb{\mathbb{N}}  
\def\MAX{{\mathcal M}}
\def\Dvv{\partial_{vv}}
\def\Dxx{\partial_{xx}}
\newcommand{\phis}{\varphi}
\def\FK{{\mathcal{F}}_{K=2}}
\newcommand{\ud}{\mathrm{d}}
\newcommand{\ui}{\mathrm{i}}
\renewcommand{\T}{\mathcal{T}}
\renewcommand{\H}{\mathcal{H}}
\renewcommand{\epsilon}{\varepsilon}

\renewcommand{\geq}{\geqslant}
\renewcommand{\leq}{\leqslant}

\newcommand{\grad}{\boldsymbol{\nabla}}

\newcommand{\cosech}{\mathop{\mathrm{cosech}}}

\newcommand{\ie}{\emph{i.e.}~}
\newcommand{\eg}{\emph{e.g.}~}

\newcommand{\half}{{\textstyle{1\over2}}}
\newcommand{\third}{{\textstyle{1\over3}}}

\title{Main topics of the NumHyp--2015' discussion session}
\date {19\up{th} June 2015}

\author[D.~Dutykh]{Denys Dutykh}
\address{CNRS -- LAMA UMR 5127, Universit\'e Savoie Mont Blanc, Camus Scientifique, 73376 Le Bourget-du-Lac, France}
\email{Denys.Dutykh@univ-savoie.fr}
\urladdr{http://www.denys-dutykh.com/\ https://dutykh.github.io/}

\author[L.~Gosse]{Laurent Gosse}
\address{CIAC-CNR, sezione di Roma, Italy}
\email{l.gosse@ba.iac.cnr.it}
\urladdr{https://sites.google.com/site/laurentgossecnr/}

\subjclass{35L65; 65M08; 76B15; 76B25}
\keywords{Conservation and balance laws; well-balanced scheme; nonlinear resonance; water waves; dispersive wave propagation}

\begin{document}
\def\smfbyname{}

\frontmatter

\begin{abstract}
Three main topics were raised in this discussion session which took place on the 19\up{th} of June at NumHyp--2015 meeting: nonlinear resonance for 1D systems of balance laws, dispersive extensions of standard hyperbolic conservation laws, and the validation of weakly dispersive shallow water wave models. An introductory overview with many bibliographic references is provided for all these topics. A numerical strategy, based on kinetic formulation, able to overcome resonance issues is presented as well as a Well-Balanced (WB) technique for Vlasov--Fokker--Planck equations are outlined. This WB scheme relies on spectral representation of stationary solutions.
\end{abstract}

\begin{altabstract}
Trois sujets principaux ont \'et\'e soulev\'es lors de la session de discussions qui a eu lieu en d\'ebut d'apr\`es-midi le 19 juin 2015 au colloque NumHyp--2015: r\'esonances nonlin\'eaires pour les syst\`emes 1D des lois de conservation, la prise en compte des effets dispersifs dans des lois de conservation et la validation des mod\`eles faiblement dispersifs pour la propagation des ondes longues. Une br\`eve r\'evu de la bibliographie est donn\'ee pour tous les sujets abord\'es. Une strat\'egie num\'erique, bas\'ee sur la formulation adapt\'ee, qui permet de contourner les r\'esonances nonlin\'eaires est pr\'esent\'ee. En plus, nous d\'ecrivons une technique bien \'equilibr\'ee pour les \'equations de Vlasov--Fokker--Planck qui est bas\'ee sur la repr\'esentation sp\'ectrale des solutions stationnaires.
\end{altabstract}

\maketitle

\newpage
\tableofcontents

\mainmatter

\section{Nonlinear resonance: issues and advances}

In the context of 1D scalar balance laws, nonlinear resonance, \ie superimposition and mutual reinforcement of convective waves with a source term of bounded extent is an old problem: see for instance, \cite{Amadori2006, Amadori2004, Seguin2012, Chen1996, Glimm1984, Goatin2004, Ha2003, Hong2006, Hong2003, Temple2004, Hua2008, Isaacson1995, Li1983, Lien1999, Liu1987}. Just to scratch the surface of this delicate subject, let us say that nonlinear resonance occurs when convective waves can see their speed of propagation vanish inside a given domain where the source term is active, possibly yielding strong amplification effects, or even general blow-up of the weak solution.

\subsection{Non-interacting asymptotic states and Temple system}

For concreteness, consider a 1D convex scalar balance law, where $g'(u)$ has no definite sign,
\begin{align}\label{eq:bal}
  \Dt u\ +\ \Dx f(u)\ &=\ k(x)g(u), \qquad t,x\ \in\ \R^+_* \times \R, \\ 
  u(t = 0,x) &= u_0(x) \in L^1 \cap BV(\R). \nonumber
\end{align}
The (possibly position-dependent) coefficient $k(x)$ is not anecdotal: indeed, it rules somehow the large-time behavior of \eqref{eq:bal} by making it switch
\begin{itemize}
  \item from a bunch of (possibly weak) traveling waves connecting (isolated) zeros of $g$ when $k(x) \equiv C \in \R^+$, a given constant, see \eg \cite{Fan1993, Harterich2005, Mascia1997},
  \item to an ordered ``\emph{scattering state}'' when $0 \leq k \in L^1(\R)$, described in \cite{Li1983, Liu1979}, consisting in a ``\emph{steady wave}'', solution of the stationary equation, close to the origin, surrounded by homogeneous, non-interacting convective waves far away from the support of $k(x)$.
\end{itemize}

Obviously, the ``\emph{steady wave}'' belonging to the scattering state can be seen as a traveling wave with zero speed. By following the canvas suggested in \cite{Glimm1986, Baiti1997}, one reformulates \eqref{eq:bal} as a peculiar $2 \times 2$ Temple system \cite{LeVeque1985}, just by introducing a fake variable $a(x)$, such that $\Dx a = k$. Requiring $k$ to be integrable means that $a \in BV(\R)$, the functions of bounded variation.
\begin{align}\label{eq:bal-sys}
  \Dt u\ +\ \Dx f(u)\ -\ g(u)\Dx a\ &=\ 0,\\ 
  \Dt a\ &=\ 0.
\end{align}
Such a reformulation is the central topic of \cite{Isaacson1995}, and more tacitly of \cite[Eqn. (3.1)]{Weinan1992}. It is the main stepping stone to build very accurate well-balanced schemes because it is homogeneous, hence given a reliable Riemann solver including the non-conservative product $g(u)\Dx a$, the Godunov scheme will result exact at numerical steady-state. Such a Godunov scheme typically contains one supplementary, very linearly degenerate, wave associated to the null eigenvalue.

According to Glimm--Sharp terminology \cite{Glimm1986}, such a reformulation of position-dependent balance laws \eqref{eq:bal}, given a source term of bounded extent, consists in handling it as a countable collection of ``\emph{local scattering centers}'' (that is to say, Dirac masses where the source's effects are localized). This meets with the usual picture of building a numerical scheme from the juxtaposition of self-similar Riemann profiles, which match the asymptotic behavior (indeed, ``\emph{a Riemann solver is Moeller's outgoing wave operator ${\mathcal W}^+$ of a scattering problem}'', \cite{Glimm1988}).

\subsection{Non-strict hyperbolicity and wave trapping effects}

A numerical scheme involving \eqref{eq:bal-sys} as its building block raises the question of its stability; for instance, starting from a smooth $k(x)$, hence a smooth $a(x)$, and its corresponding piecewise-constant approximation $P^\DX a$, given a mesh-size $\DX > 0$, what happens to a sequence of solutions $u^\epsilon$, $a^\epsilon$, where $a^\epsilon \to P^\DX a$ as $\epsilon \to 0$? This question (of continuous dependence with respect to data) was eluded in \cite{Isaacson1995, Weinan1992, Greenberg1996} but fully addressed in \cite[\S 3.3]{Gosse2001}, relying on previous estimates in \cite{Mascia1999} and convenient $BV$ bounds.

Compactness holds under the ``\emph{non-resonance assumption}'', roughly $f'(u) \neq 0$ (see also \cite{Natalini1994} or even \cite{Guerra2004}), essentially for two reasons:
\begin{enumerate}
  \item In order to take full advantage of Temple's theory, \ie Riemann invariants seeing their total variation decay in time \cite{LeVeque1985}, one needs the mapping $(u,a) \mapsto (w,a)$, with $w$ the Riemann invariant written in \cite[Eqn. (37)]{Gosse2001}, to be a diffeomorphism.
  \item Later, to exploit \cite{Mascia1999} and give a sense to the product $g(u)\Dx\, a$ when $a$ is discontinuous, one needs ``\emph{outgoing traces}'' which are available if characteristic velocities don't vanish.
\end{enumerate}
Worse, it was found that $L^1$ stability estimates, with respect to both time and data, blow up when resonance occurs because some constants are of the order of $|1/f'|$: see \cite{Amadori2004, Amadori2013, Guerra2004}.

For 1D nonlinear systems of conservation laws, both Glimm's theory and $L^1$ stability mostly rely on strict hyperbolicity, which constraints the corresponding dynamics to decay, as times passes, from a very complex picture of interacting elementary waves, to a well-ordered, non-interacting, scattering state involving only non-approaching waves. Such a decay takes place as Glimm's interaction potential slowly decays, and both information and entropy are dissipated in shocks \cite{Liu1979, Glimm1984}. Instead, nonlinear resonance corresponds to a formal merging of both the static and convective characteristic families in \eqref{eq:bal-sys}, namely $f'(u) = 0$. The classical Lax's theorem doesn't ensure existence and uniqueness for the Riemann problem. In the scalar case, despite non-uniqueness issues, a quite complete catalog of Riemann solutions was published in \cite{Isaacson1995}, together with numerical applications in \cite{Noussair2000}. When it comes to general 1D systems, the picture worsens: see \cite{Goatin2004, Hong2006, Hong2003, Temple2004, Hua2008}. Stability was further studied in \cite{Lien1999, Ha2003}.

In \cite[Ex.~11, page~71]{Bressan2009}, Bressan presents an illuminating counter-example which reveals why one should not expect any continuous dependence of $BV$ solutions to a simple system of the form \eqref{eq:bal-sys} when strict hyperbolicity fails (in the sense that eigenvalues coalesce). Accordingly,
\begin{align*}
  \Dt u\ +\ \Dx f(a,u)\ &=\ 0, \quad f(a,u)\ =\ 8au - u^2, \\
  \Dt a\ &=\ 0,
\end{align*}
models car traffic on a highway, with $a(x)$ the number of lanes opened. Eigenvalues are $\Lambda = \{0, 8a - 2u\}$, so resonance occurs when $u = 4a$. If some repair works take place, $a(x)$ switches from $2$ to $1$: it is thus possible to cook up an admissible Riemann fan, containing a resonant state $a = 1$, $u = 4$ and continuous dependence with respect to initial data is violated.

A more favorable situation emerges when only a linearly degenerate field superimposes with the standing wave: this is essentially what is discussed in \cite{Desveaux2015}. Ill-posedness seems to be confined to resonance between genuinely nonlinear fields and the standing wave. Besides, in the scalar case, it appears that an astute change of variables, originally proposed in \cite{Liu1987}, allows to set up a resonant wave front tracking algorithm \cite{Karlsen2004}, or a Godunov scheme \cite{Amadori2004, Isaacson1995, Noussair2000}.

\subsection{A numerical case-study for isentropic Euler with $\gamma = 3$}

Hereafter, we follow \cite[Ch. 6 pp.~107-110]{Gosse2013} (see also \cite{Perthame2001}) in order to take fully advantage of a particular kinetic formulation for a genuinely nonlinear, non-Temple, hyperbolic system of balance laws. Let the initial-boundary value problem (IBVP) be posed in $x \in (-1, 1)$, $t > 0$,
\begin{align}\label{eq:EP1}
  \Dt\rho\ +\ \Dx (\rho u) \ &=\ 0, \\ 
  \Dt(\rho u)\ +\ \Dx\bigl(\rho u^2 + \frac{\rho^3}{12}\bigr)\ +\ \rho \Dx \phi\ +\ \frac{\rho u}{\tau(x)}\ &=\ 0, \label{eq:EP2}
\end{align}
where the electric field $E = -\Dx \phi$ and $\phi$ satisfies a repulsive Poisson equation:
\begin{align}\label{eq:PE}
  \lambda(x)\Dxx\phi\ &=\ \rho_D(x)\ -\ \rho, \\
  \phi(t,\; x = -1)\ &=\ 0, \nonumber \\
  \phi(t,\; x = 1)\ &=\ -V \geq 0. \nonumber
\end{align}
The positive parameters $\tau(x)$, $\lambda(x)$ are the damping coefficient and scaled Debye's length: they are related to the material inside which electronic conduction takes place. Generically, a strongly doped material induces both a small $\tau$ and $\lambda$. The doping profile is given by $\rho_D$: a rescaling of the system allows to get $\rho_D \in (0,1)$. Thus, in a $n^+ n n^+$ device, we should always get discontinuous parameters like (for instance) in \cite[page 235]{Gosse2013}. The pressure law with adiabatic exponent $\gamma = 3$ has the nice property that the quasi-linear system \eqref{eq:EP1}, \eqref{eq:EP2} admits a ``\emph{pure}'' kinetic formulation, called ``\emph{$K$-multi-branch entropy solution}'' after \cite{Brenier1998}. Here, we just pick $K = 2$ in order to realize the moments $\rho, \rho u$, so the kinetic density $f(t,x,v)$ solves,
\begin{equation*}
  \Dt f\ +\ v\cdot\Dx f\ +\ E\Dv f\ =\ -\Dvv \mu, \quad 
  -\Dvv \mu\ =\ \lim_{\epsilon \to 0}\frac{1}{\epsilon}\Big(\MAX_{K = 2}(v;\; \rho,\, \rho u)\ -\ f\Big),
\end{equation*}
and the ``\emph{2-branch Maxwellian}'' is given by two Heaviside functions, (see first pages of \cite{Brenier1998})
\begin{equation}\label{eq:maxell}
  \MAX_{K = 2}(v;\; \rho,\, \rho u)\ =\ Y(v-u^-)\ -\ Y(v-u^+), \qquad u^\pm\ =\ u\ \pm\ \frac{\rho}{2},
\end{equation}
where $u^\pm$ are Riemann invariants for \eqref{eq:EP1}, \eqref{eq:EP2}: 
\begin{equation*}
  \rho\ =\ u^+ - u^-, \qquad
  \rho u\ =\ \frac{1}{2}\bigl(|u^+|^2 - |u^-|^2\bigr).
\end{equation*}
This property was extensively used in \cite{Vasseur1999} to prove regularity properties for $L^p$ weak solutions of the homogeneous isentropic Euler system with $\gamma = 3$. Here, the main goal with \eqref{eq:EP1}--\eqref{eq:PE} is to derive the so--called ``\emph{current-voltage}'' relation for the device under consideration: that is to say, a curve giving the stationary circulating current as a function of the potential drop,
\begin{equation*}
  V \mapsto J_{\mathrm{stat}}(x), \quad \mbox{ with }\quad J_{\mathrm{stat}}(x) = \rho u (x) \mbox{ is the (constant) stationary current}. 
\end{equation*}
Hence it is highly desirable to have at hand a method able to deliver a constant momentum at numerical steady-state: this is a similar situation to a ``\emph{shallow water system with a $\rho$-dependent topography}''. However, such methods do not really exist, see for instance \cite{Ballestra2004}.

\subsection{A fictitious (but resonant) model of 1D semiconductor}

We intend to build the numerical method on the following approximate kinetic formulation,
\begin{align}\label{eq:vla}
  \Dt f + v\cdot\Dx f + E(t,x) \Dv f\ &=\ \sum_{n \in \Nbb} \Big(\MAX_{K = 2}(v;\rho,\rho u)-f\Big)\delta(t-n\DT), \\ 
  E\ &=\ -\Dx\phi, \nonumber
\end{align}
first when $\tau \to + \infty$ (no damping), which reduces to a time-marching, splitting algorithm:
\begin{itemize}
  \item In the open layers $(t^n, t^{n + 1})$, a collision-less Vlasov equation is solved,
  \item At discrete times $t^n = n\DT$, $f$ is projected onto its local 2-branch Maxwellian \eqref{eq:maxell}.
  \item Finally, the self-consistent Poisson equation is solved by standard finite-differences.
\end{itemize}
The resulting well-balanced numerical fluxes are given for any $K \in \Nbb$ in \cite[(6.31) p.~110]{Gosse2013}.
\begin{align*}
  \rho(t,x)\ &=\ \int_\R f(t,x,v)\;\ud v\ =\ (u^+ - u^-)(t,x), \\ 
  u(t,x)\ &=\ \frac{\int_\R v f(t,x,v)\;\ud v}{\int_\R f(t,x,v)\;\ud v}\ =\ \frac{(u^+ + u^-)(t,x)}{2}.
\end{align*}
Let $f(\rho,\rho u)$ (or equivalently $f(u^+,u^-)$) stand for the flux function in \eqref{eq:EP1}, \eqref{eq:EP2}:
\begin{equation*}
  f(\rho,\,\rho u)\ =\ \int_\R (v, v^2)\MAX_{K=2}(v;\; \rho,\, \rho u)\;\ud v = \left(\frac{|u^+|^2-|u^-|^2}{2},\; \frac{(u^+)^3-(u^-)^3}{3} \right).
\end{equation*}
A kinetic flux-vector splitting scheme rests on the definition of $f^\pm$ as follows,
\begin{equation*}
  f^\pm(\rho,\rho u)\ =\ \left(\int_0^{\pm \max(0,\pm u^+)} - \int_0^{\pm \max(0,\pm u^-)}\right)  (v, v^2)\MAX_{K=2}(v;\;\rho,\, \rho u)\;\ud v
\end{equation*}
Letting $\vec m_j^n = (\rho_j^n, \rho u_j^n)$ stand for the numerical approximation of the moment vector,
\begin{equation}\label{PS-K-scheme}
  \vec m_j^{n+1}\ =\ \vec m_j^n\ -\ \frac{\DT}{\DX}\left[(\FK^-)_{j+\frac 1 2} - f^-(\vec m_j^n) + f^+(\vec m_j^n) - (\FK^+)_{j-\frac 1 2}\right],
\end{equation}
with $\vec u = (u^+,u^-)$, $\Delta \phi$ be Riemann invariants and potential jump at each interface,
\begin{equation}\label{eq:PS-flyx}
\left\{\begin{array}{l}
  f^\pm(\vec m)\ =\ f(\pm \max(0,\pm \vec u)),\quad \to \quad \mbox{there is a typo in \cite[(6.31)]{Gosse2013},}\\
  \begin{array}{rcl}
    \FK^+&=&f \left(\sqrt{\max\bigl(0,(\max(0,\vec u_{\mathrm{left}})^2-2\Delta \phi\bigr)}\right)\\
    &&\qquad -f \left(\min(\max(0,-\vec u_{\mathrm{right}}),\sqrt{\max(0,-2\Delta \phi)})\right),
  \end{array}\\
  \begin{array}{rcl}
    \FK^-&=&f \left(-\sqrt{\max\bigl(0,(\min(0,\vec u_{\mathrm{right}})^2+2\Delta \phi\bigr)}\right)\\
    &&\qquad -f \left(-\min\bigl(\max(0,\vec u_{\mathrm{left}}),\sqrt{\max(0,2\Delta \phi)}\bigr)\right).
  \end{array}\\
  \end{array}\right.
\end{equation}
In \eqref{eq:PS-flyx}, each $\FK^\pm$ contains 2 terms: the first one is the standard well-balanced flux, where the effects of the potential jump $\Delta \phi$ is included in the convective terms. This first term can be derived with standard considerations without invoking any kinetic formulation. However, there's now a second term, specific to the kinetic formalism, which accounts for particles whose kinetic energy is too low to pass through the potential barrier; hence they are reflected (like in \cite[p.~210]{Perthame2001}) and this term strongly contributes when nonlinear resonance (that is, sonic points) occurs within numerical simulations of the tricky model \eqref{eq:EP1}--\eqref{eq:PE}.

\begin{itemize}

\item the kinetic formulation involving the Vlasov model \eqref{eq:vla} leads to stationary regimes usually endowed with 2 sonic points: one in a rarefying region located inside the channel, the other inside the drain, but attached to a strong transonic shock: see Fig.~\ref{trans}. Despite the well-balanced character of the scheme, there is a strong spike in the stationary current located at the place of the transonic shock (the Mach number $\frac{2u}{\rho}$ passes from $8.5$ to $0.2$), because of shock fitting on the Cartesian grid (which is coarse, with $2^6$ points in the interval $(-1,1)$).

\begin{figure}
  \centering
  \includegraphics[width=0.48\linewidth]{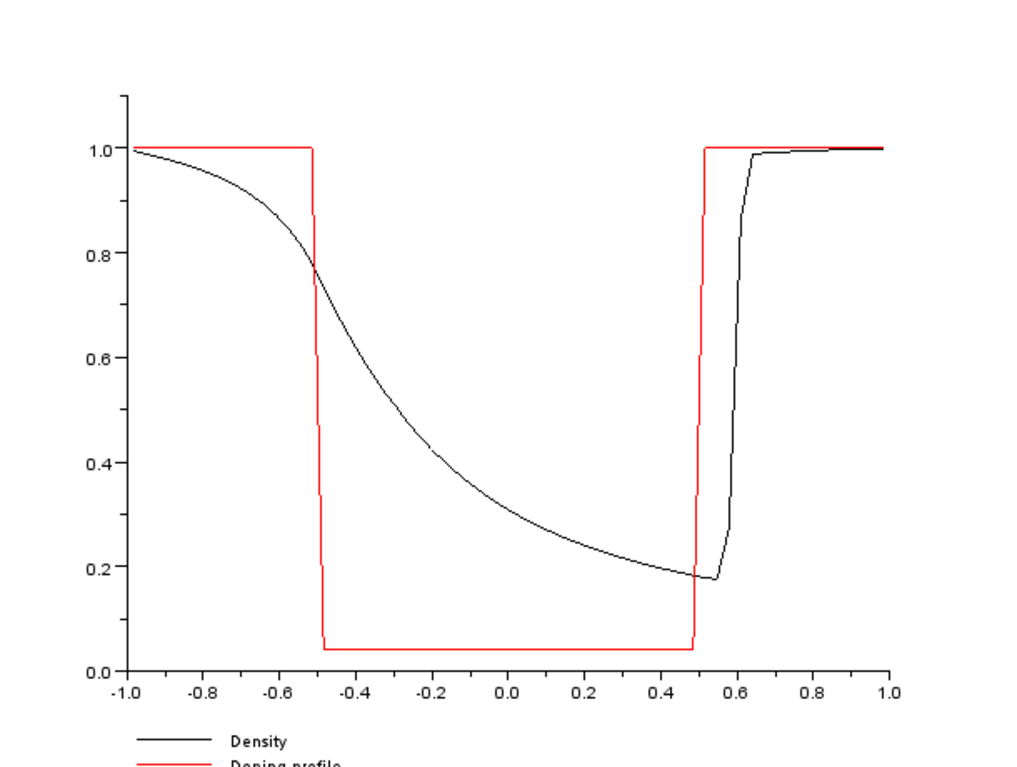}
  \includegraphics[width=0.48\linewidth]{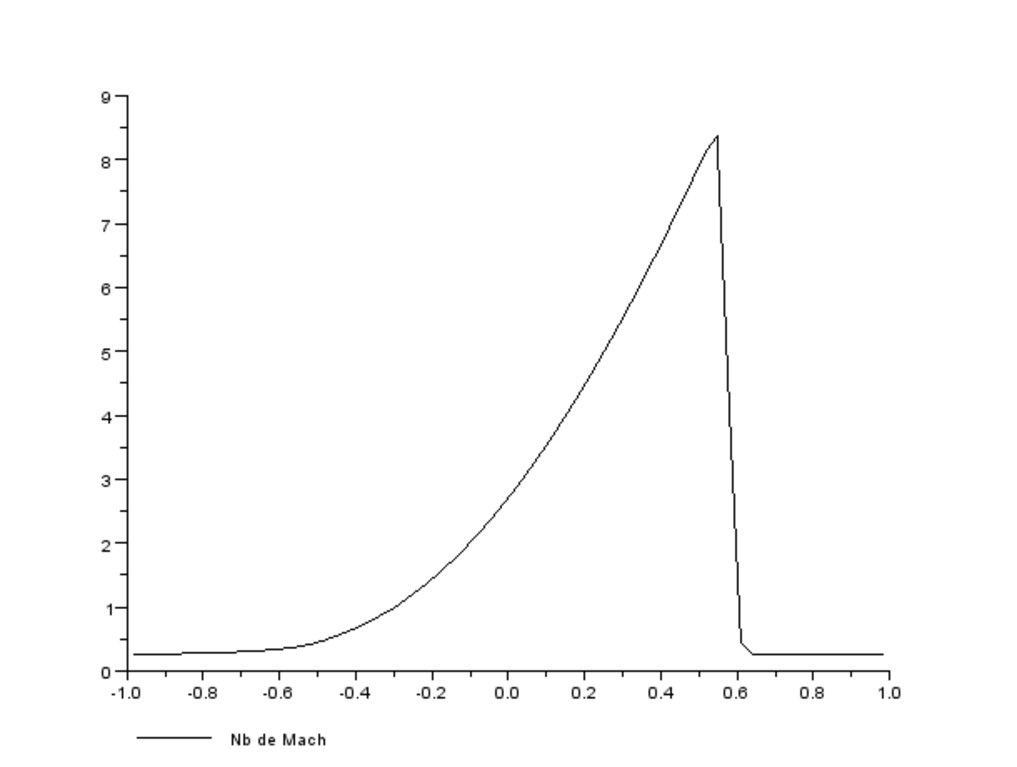}
  \includegraphics[width=0.48\linewidth]{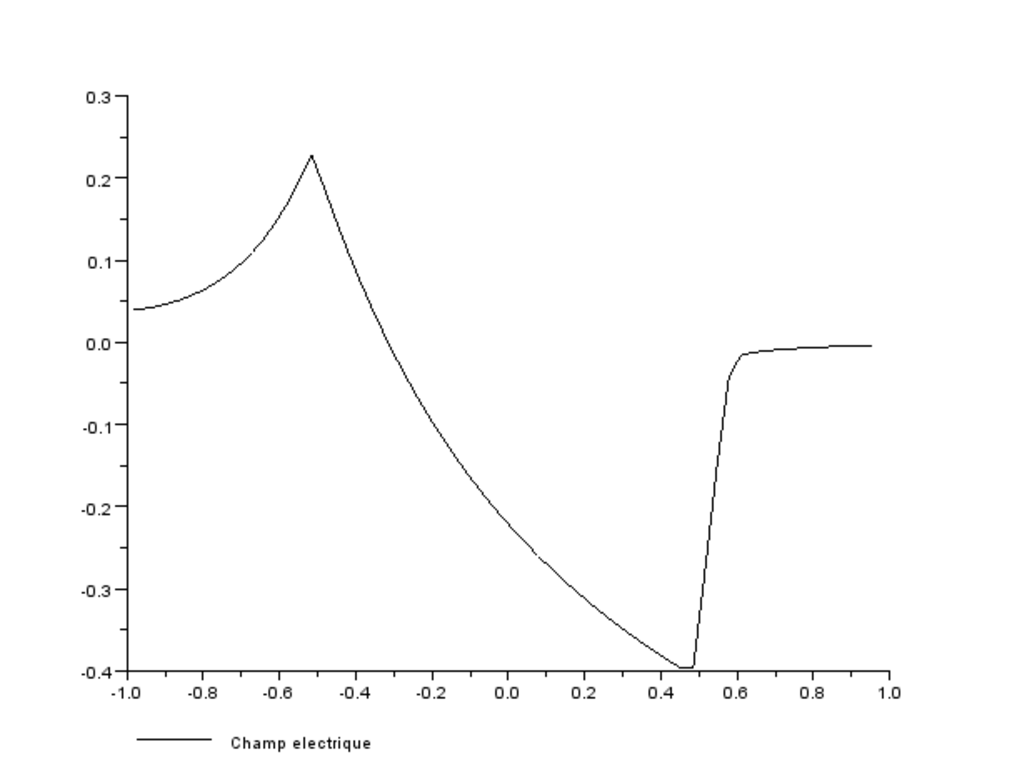}
  \includegraphics[width=0.48\linewidth]{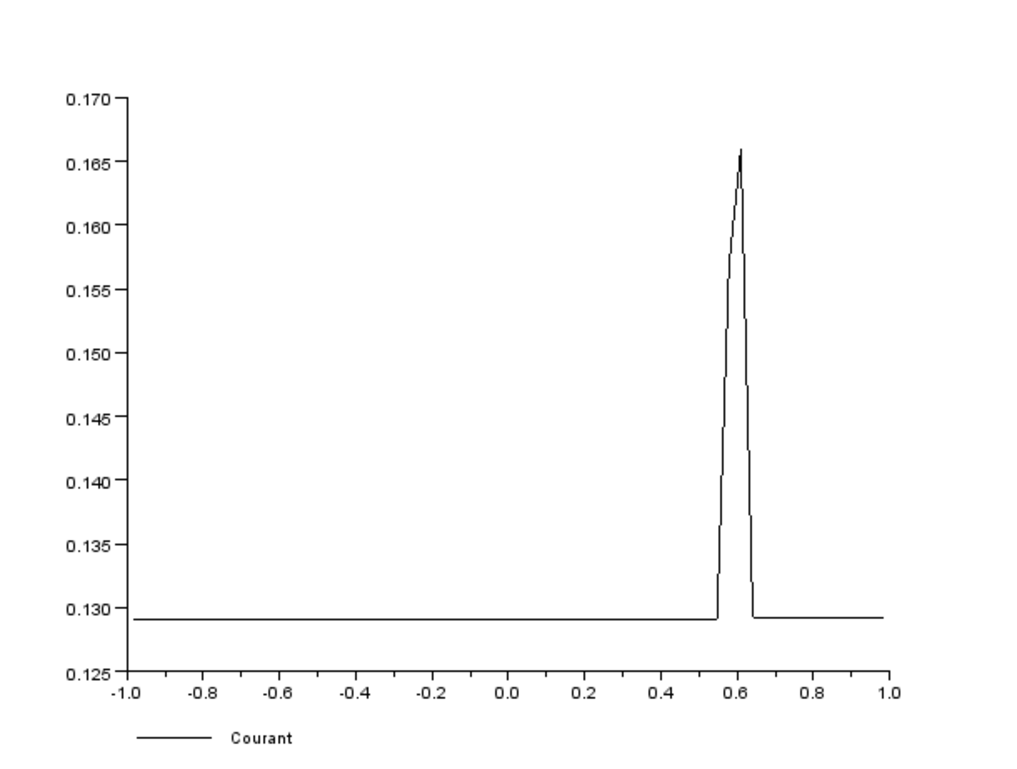}
  \caption{\small\em Numerical results for $V = -0.15$, $\lambda = 0.15, 0.5$ and no damping.}
  \label{trans} 
\end{figure}

\item This feature shows the shortcomings of ``\emph{lake at rest}'' ideas, since imposing zero volt at the edges of the device $V = 0$ brings a perfect numerical (smooth) solution: see Fig.~\ref{zero}.

\begin{figure}
  \centering
  \includegraphics[width=0.48\linewidth]{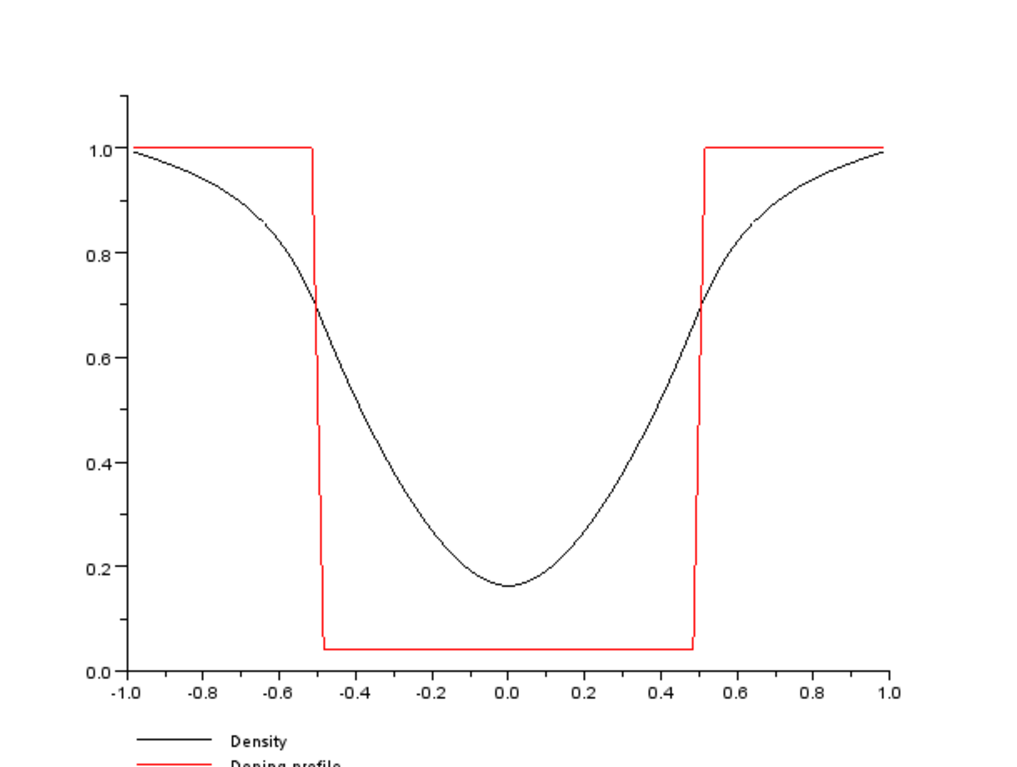}
  \includegraphics[width=0.48\linewidth]{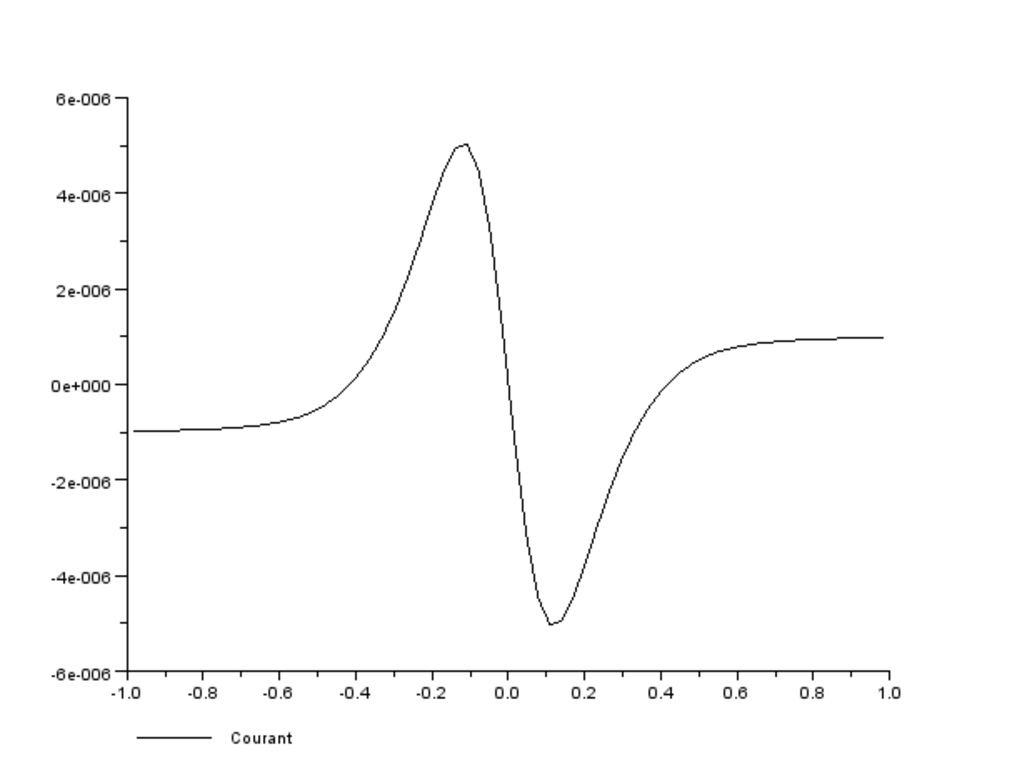}
  \caption{\small\em Numerical results on $\rho, \rho u$ for $V=0$, $\lambda=0.15, 0.5$ and no damping.}
  \label{zero}
\end{figure}

\begin{figure}
  \centering
  \includegraphics[width=0.48\linewidth]{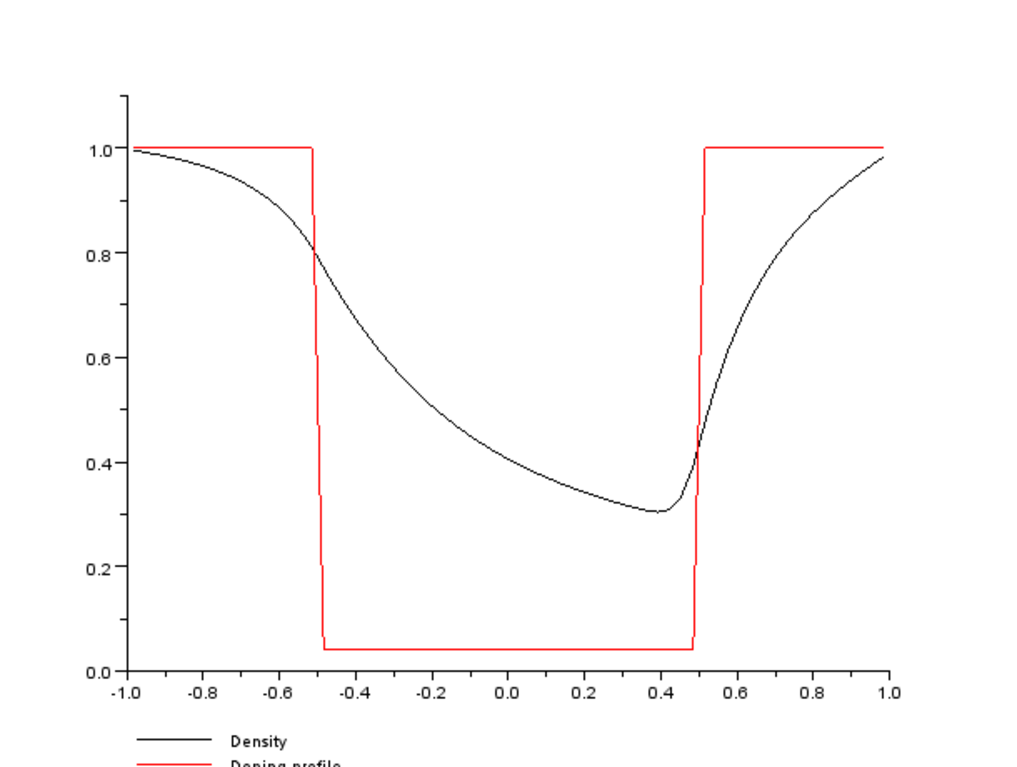}
  \includegraphics[width=0.48\linewidth]{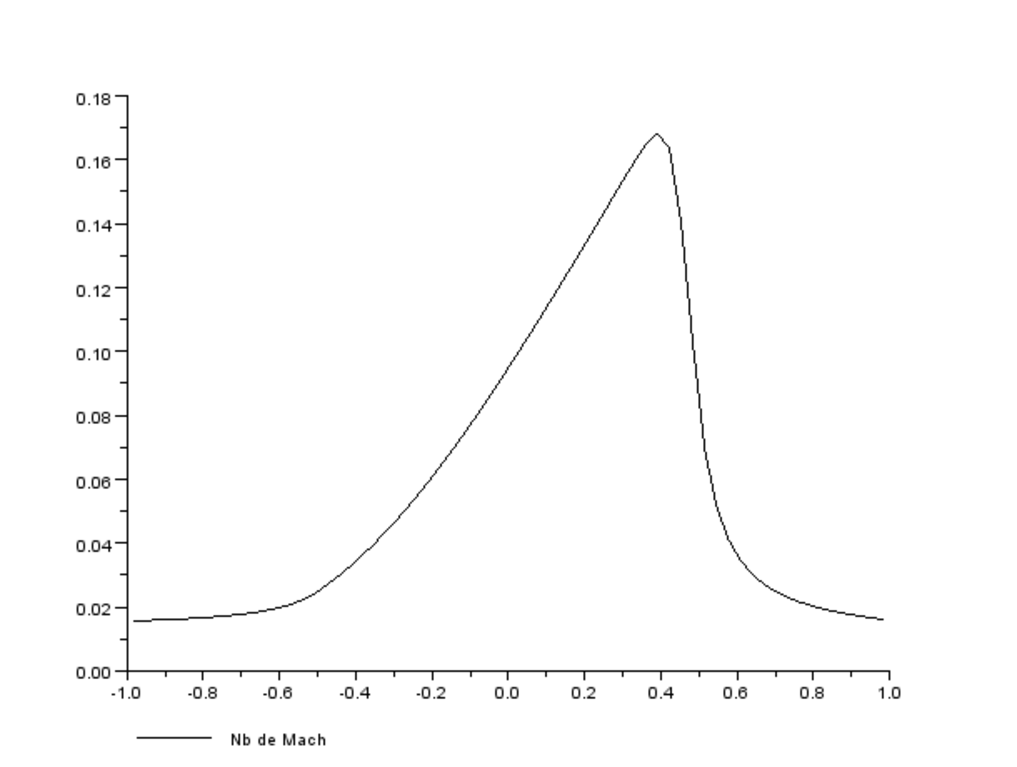}
  \includegraphics[width=0.48\linewidth]{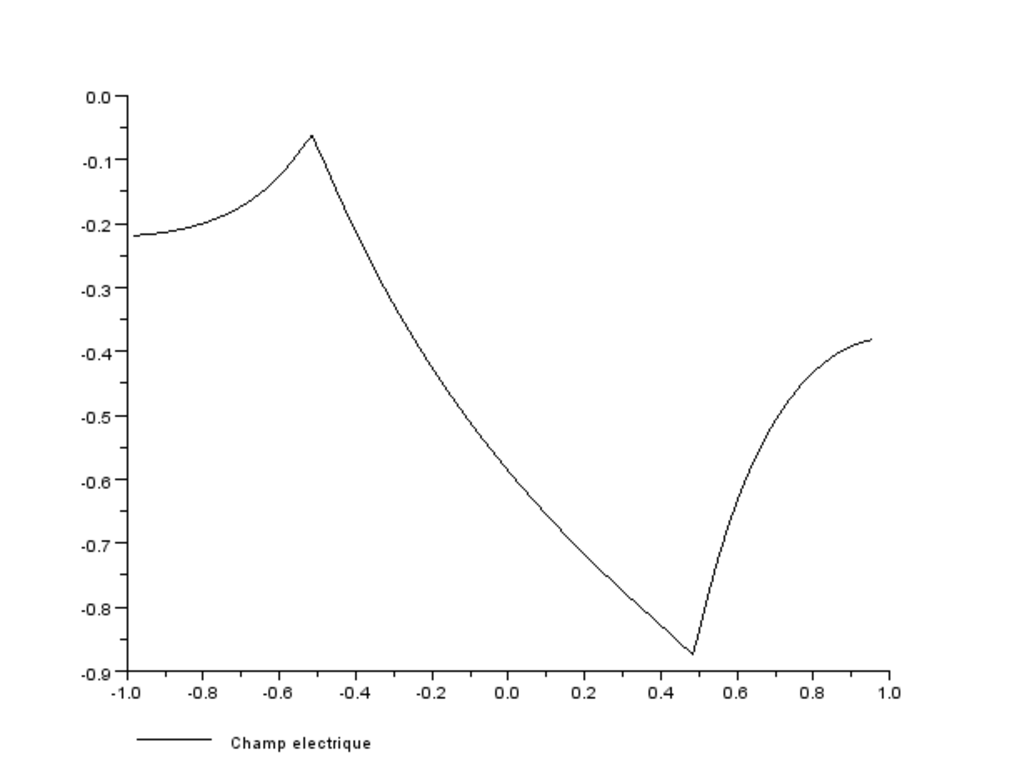}
  \includegraphics[width=0.48\linewidth]{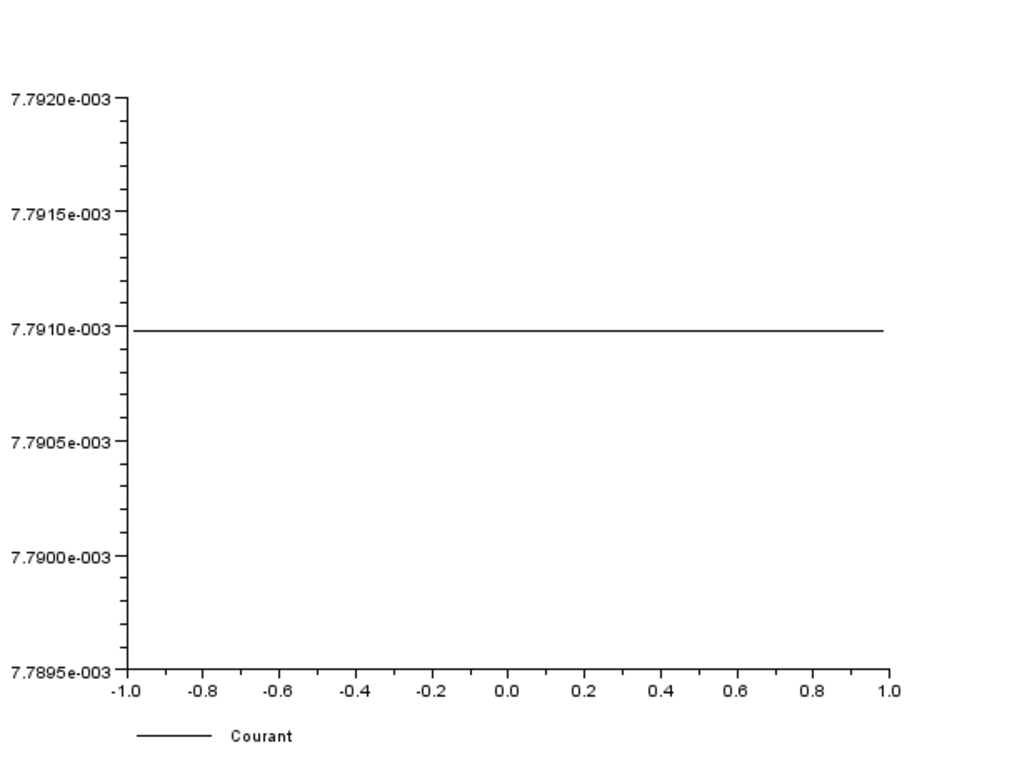}
  \caption{\small\em Numerical results for $V=-0.9$, $\lambda=0.15, 0.5$ and damping $\tau \equiv 1$.}
  \label{damp}
\end{figure}

\item Now, the damping term $\frac{\rho u}{\tau}$ is helpful in order to tame the dynamics, as it cancels sonic shocks for moderate values of $V$. A simple manner to include the damping effect into the kinetic scheme consists in substituting the ``\emph{potential jumps}'' in \eqref{eq:PS-flyx} with ``\emph{augmented potential jumps}'' which contain a supplementary term $-\frac{1}{2\tau} (u_{\mathrm{left}} + u_{\mathrm{right}})$. However, this cannot handle discontinuous dampings $\tau(x)$. For uniform, constant dampings, a perfect stabilization emerges, without sonic shocks: see Fig.~\ref{damp} where $\tau \equiv 1$ and $V = -0.9$.
\end{itemize}

Despite its simplicity, such an example can be extended to more involved moment systems based on the formalism of $K$-multibranch solutions, but one quickly needs the so--called ``\emph{Markov moment inversion}'', see \cite[page 1621]{Gosse2008}, in order to express numerical fluxes. Results in \cite{Lien1999} indicate that one may expect stability in resonant cases only when the source (a ``\emph{sink}'') dissipates strongly in order to prevent trapped waves from being amplified beyond control.

\section{Disperse flows and hyperbolic-kinetic couplings}

Transport of pollutants in a river can be modeled, in first approximation, by means of a supplementary conservation law which induces a contact discontinuity moving at the material velocity in the Riemann problem. However, the characteristic scale of the considered particles suggests that a hyperbolic-kinetic coupling is probably more adequate, as it is of current use in several spray modelings and simulations. We follow \cite{Caflisch1983, Carrillo2006, Goudon2014, Patankar2001, Patankar2001a} and \cite[Ch.~12]{Gosse2013}.

\subsection{Sediment transport in shallow water models}

An elementary 1D coupled system of pollutant transport \cite{Heemink1990} (see also \cite{Audusse2003, Chertock2006}) is described,
\begin{equation*}\label{sha-pol}
  \begin{array}{rcl}
  \Dt\rho\ +\ \Dx(\rho u) &\ =\ & 0, \\ 
  \Dt(\rho u)\ +\ \Dx (\rho u^2 +\frac{\rho^2}{2} )&\ =\ & \int_\R (v-u)f(t,x,v)\;\ud v, \smallskip \\
  \Dt f\ +\ v \cdot \Dx f &\ =\ & \Dv \big[(v-u)f - \Dv f\big],
  \end{array}
\end{equation*}
inspired by \eg \cite{Goudon2014}, and where gravity potentials or various dimensional parameters were ignored for simplicity. The kinetic density $f(t,x,v)$ stands for the disperse, rarefied phase, the pollutant, explicitly depending on a microscopic velocity variable, $v \in \R$. Beside the usual Fokker--Planck friction term, the right-hand side of the kinetic equation displays an acceleration term involving the macroscopic velocity $u(t,x)$, so that pollutants evolve according to both a drift effect of the underlying flow and a stochastic diffusion in $v$. Both macroscopic and microscopic descriptions are weakly coupled by means of their right-hand sides which renders a drag effect of one phase onto the other. Within a standard spitting approach, Fokker--Planck and shallow water models are solved independently of each other, so they inevitably involve position-dependent coefficients inside their respective source terms,
\begin{align*}
  u(t^n,x) &\mbox{ in the kinetic equation}, \\ 
  f(t^n,x,\cdot) &\mbox{ in the momentum equation},
\end{align*}
indicating that a well-balanced strategy may be desirable.

\subsection{A simple Burgers/Vlasov--Fokker--Planck toy model}

A more tractable model is obtained by reducing the shallow water system to Burgers equation, 
\begin{align*}
  \Dt u\ +\ \Dx\bigl(\frac{u^2}{2}\bigr)\ &=\ \int_\R (v-u)f(t,x,v)dv, \\
  \Dt f\ +\ v \cdot \Dx f\ &=\ \Dv \big((v-u)f - \kappa \Dv f\big),
\end{align*}
as suggested in \cite{Domelevo2001}, where existence and stability results were established. Such a weakly coupled system is simple enough to be numerically solved within an overall well-balanced approach, itself based on a Godunov discretization involving the reformulation \eqref{eq:bal-sys} for the macroscopic part, and the so-called ``\emph{elementary solutions}'' framework for the kinetic one.

It is well-known that a good knowledge of stationary ``\emph{standing wave}'' solutions is a decisive advantage when setting up WB discretizations, a way to generate extra-accuracy \cite{Amadori2013, Amadori2015} when position-dependent coefficients are involved. By rewriting one right-hand side as,
\begin{align*}
  \int_\R (v-u)f(t,x,v)dv\ &=\ J(t,x)\ -\ u(t,x)\int_\R f(t,x,v)\;\ud v, \\ 
  J(t,x)\ &=\ \int_\R vf(t,x,v)\;\ud v,
\end{align*}
given $f$, it is possible to solve the stationary equation by means of Lambert's $W$-function. Yet, the stationary solutions of various linear kinetic equations (with incoming boundary conditions) was the object of intensive research (sometimes called ``\emph{Caseology}'', after the seminal 1960 paper by Kenneth \textsc{Case}) for several decades: see \eg the books \cite{Case1967, Cercignani1969, Greenberg1984, Kaper1982}.

In contrast with kinetic equations rendering neutron transport or linearized BGK models \cite[Ch. VII]{Cercignani1969}, the stationary problem for Fokker--Planck equations in slab geometry can be solved analytically by means of Sturm--Liouville theory and Hermite functions, \cite{Beals1979, Beals1983, Burschka1982, Cercignani1992, Pagani1970}. Moreover, the inclusion of a constant force field (here $u(t,x)$, approximated in each cell by a piecewise constant function) doesn't ask for big changes in spectral representations; this feature is specific to Fokker--Planck models as the ones related to integral collision terms hardly support forcing terms expressed by first-order derivatives in $v$, see \cite{Dalitz1997, Toepffer1997}.

In the realm of a discrete-ordinates approximation in the $v$ variable, let us consider 
\begin{equation}\label{eq:VFP-FP-vla-eqn}
  \Dt f + v\cdot \Dx f + u\cdot \Dv f\ =\ \Dv (f+\kappa \Dv f),\qquad u \in \R, \qquad \kappa \in \R_*^+.
\end{equation}

By assuming a rather standard translation-invariant ansatz $\psi(x,v) = \exp(-\lambda x -\mu v)\varphi(v)$ in the context of the stationary problem for \eqref{eq:VFP-FP-vla-eqn}, one is led to a Sturm--Liouville problem,
\begin{equation*}
  \left(\kappa\frac{\ud^2}{\ud v^2} + (v - 2\mu\kappa - u)\frac{\ud}{\ud v} + (\lambda - \mu)v + \mu^2 \kappa+ \mu u\right)\varphi(v)\ =\ 0.
\end{equation*}
The rescaled Sturm--Liouville operator which is used in the absence of any forcing field is:
\begin{equation*}
  S_\kappa(v)\ :=\ \frac{\ud}{\ud v}\left(v(\cdot) + \kappa \frac{\ud(\cdot)}{\ud v}\right)\ =\ \kappa \frac{\ud^2 (\cdot)}{\ud v^2}\ +\ v\frac{\ud(\cdot)}{\ud v} + (\cdot)
\end{equation*}
In order to recover an expression involving $S_\kappa$, the natural choice is: 
\begin{equation*}
 \lambda\ =\ \mu, \quad \kappa \mu^2 + u\mu -n\ =\ 0\quad \mbox{ for }\quad n \in \Nbb.
\end{equation*}
It yields $S_\kappa(v-2\mu \kappa-u)[\varphi] = 0$, so the eigenvalues read for any $n \in \Nbb$:
\begin{equation*}
  \mu_{\pm n}\ =\ \frac{-u \pm \sqrt{u^2 + 4\kappa n}}{2\kappa} \stackrel{u \to 0}{\longrightarrow}\ \pm\ \sqrt{\frac n \kappa}.
\end{equation*}
The eigenfunctions for \eqref{eq:VFP-FP-vla-eqn}, originally published in \cite{Burschka1982}, can be expressed as:
\begin{equation*}
  \Psi_{\pm n}(x,v)\ =\ \exp\left(-\mu_{\pm n}[x+ v]\right)H_n\left(\tilde v_{\pm n}\right)\exp\left(-\tilde v_{\pm n}^2\right),
\end{equation*}
being $H_n$ the $n^{th}$ Hermite polynomial and the translated velocities given by,
\begin{equation*}
  \tilde v_{\pm n}\ =\ \frac{v-2\mu_{\pm n} \kappa-u}{\sqrt{2\kappa}}\ =\ \frac{v \mp\sqrt{u^2 + 4\kappa n}}{\sqrt{2\kappa}}\ \stackrel{u \to 0}{\longrightarrow}\ \frac{v}{\sqrt{2\kappa}}\ \mp\ \sqrt{2n}.
\end{equation*}
It remains to derive both the two Chapman--Enskog eigenfunctions associated to $n = 0$. However, even for $u \neq 0$, they are easily found because plugging $n = 0$ leads to:
\begin{equation*}
\left\{
\begin{array}{rcl}
  \mu_0^+& = &\displaystyle \frac{-u+|u|}{2\kappa}=-\frac{ \inf(u,0)}{\kappa}, \qquad \tilde v^+_0\ =\ \frac{v-|u|}{\sqrt{2\kappa}}, \\
  \mu_0^-& = &\displaystyle \frac{-u-|u|}{2\kappa}=-\frac{ \sup(u,0)}{\kappa}, \qquad \tilde v^-_0\ =\ \frac{v+|u|}{\sqrt{2\kappa}}.
\end{array}\right.
\end{equation*}
Thus the ``\emph{diffusion eigenfunctions}'' corresponding to $n = 0$ and $u \neq 0$ read:
\begin{equation*}
\begin{array}{rcl}
  \Psi_0^+(x,v)&=&\exp\left(\inf(u,0)\frac{(x+ v)} \kappa -\frac{(v-|u|)^2}{2\kappa} \right)\\
  &=&\exp\left( -\frac{(v-|u|)^2}{2\kappa} \right)\chi_{E>0}\ +\ \exp\left(-\frac{u^2}{2\kappa}\right)\exp\left(\frac{ux} \kappa -\frac{v^2}{2\kappa} \right)\chi_{E<0},\\
  \Psi_0^-(x,v)&=&\exp\left(\sup(u,0)\frac{(x+ v)} \kappa -\frac{(v+|u|)^2}{2\kappa} \right) \quad \to \mbox{ there is a typo in \cite[p.251]{Gosse2013}},\\
  &=&\exp\left( -\frac{(v+|u|)^2}{2\kappa} \right)\chi_{u<0}\ +\ \exp\left(-\frac{u^2}{2\kappa}\right)\exp\left(\frac{ux} \kappa -\frac{v^2}{2\kappa} \right)\chi_{E>0}.
\end{array}
\end{equation*}
Finally, any stationary solution $\bar f(x,v)$ of \eqref{eq:VFP-FP-vla-eqn} spectrally decomposes into,
\begin{equation*}
  \bar f(x,v)\ =\ \underbrace{\alpha \Psi_0^+(x,v) + \beta \Psi_0^-(x,v)}_{\mbox{macroscopic behavior}}\ +\ \underbrace{\sum_{n \in \Nbb^*} A_n \Psi_n(x,v)+B_n\Psi_{-n}(x,v)}_{\not = 0 \mbox{ in the Knudsen layers}}.
\end{equation*}
With these expressions at hand, it suffices to follow the calculations in \cite[Ch. 12]{Gosse2013} to derive the scattering matrices which lead to a kinetic well-balanced scheme for \eqref{eq:VFP-FP-vla-eqn}; an example where this technique was set up in a delicate context of high-field limits is provided in \cite{Gosse2015}.

\section{Dispersive wave modelling}

In this part we shall focus on free surface wave modelling as a practically important example of nonlinear dispersive wave propagation. As it was pointed out by R.~Feynman, the usual water waves are \emph{in no respect like} the light (electro-magnetics) or acoustics (pressure waves) \cite{Feynman2005}.

\subsection{Water wave propblem}

We consider surface gravity waves in an ideal incompressible and irrotational fluid of finite depth. A two-dimensional Cartesian coordinate system $xOy$ is chosen such that the ordinate axis $Oy$ is directed vertically upwards and the horizontal axis $Ox$ coincides with the still water level. The fluid layer is bounded above by the free surface $y = \eta(x,t)$ and below by the flat bottom $y = -h$. The governing equations for the flow are \cite{Stoker1957}
\begin{eqnarray}
  \Delta\phi = \phi_{xx} + \phi_{yy} &=& 0, \qquad -h < y < \eta(x,t), \label{eq:lapl} \\
  \eta_t + \phi_x\cdot\eta_x &=& \phi_y, \qquad y = \eta(x,t), \label{eq:kin} \\
  \phi_t + \half|\grad\phi|^2 + gy &=& 0, \qquad y = \eta(x,t), \label{eq:dyn} \\
  \phi_y &=& 0, \qquad y = -h, \label{eq:bot}
\end{eqnarray}
where $\phi(x,y,t)$ is the velocity potential and $g$ is the acceleration due to gravity. The Laplace equation \eqref{eq:lapl} expresses the combination of fluid incompressibility and flow irrotationality.  This equation is completed by the boundary conditions \eqref{eq:kin}--\eqref{eq:bot}. There is one kinematic \eqref{eq:kin} and one dynamic isobarity condition \eqref{eq:dyn} on the free surface. On the solid bottom we require that the impermeability condition \eqref{eq:bot} is satisfied.

The set of equations \eqref{eq:lapl}--\eqref{eq:bot} is referred to as the full water wave problem. It has been shown to provide an excellent description for free surface flows \cite{Stoker1957}, despite its apparent simplicity.

\subsection{The method of holomorphic variables}

One of the main difficulties of the water wave problem \eqref{eq:lapl} -- \eqref{eq:bot} is that the fluid domain is unknown a priori and has to be determined amongst other unknowns. Consequently, in order to simplify the solution procedure, it would be advantageous to transform the dynamic physical domain into a fixed computational one. For numerical purposes the idea to use time-dependent conformal maps was formalized and implemented for the first time by A.~\textsc{Dyachenko} \emph{et al.} (1996) \cite{Dyachenko1996}. However, this formulation was put forward for the first time presumably by L.~\textsc{Ovsyannikov} (1974) \cite{Ovsyannikov1974}.

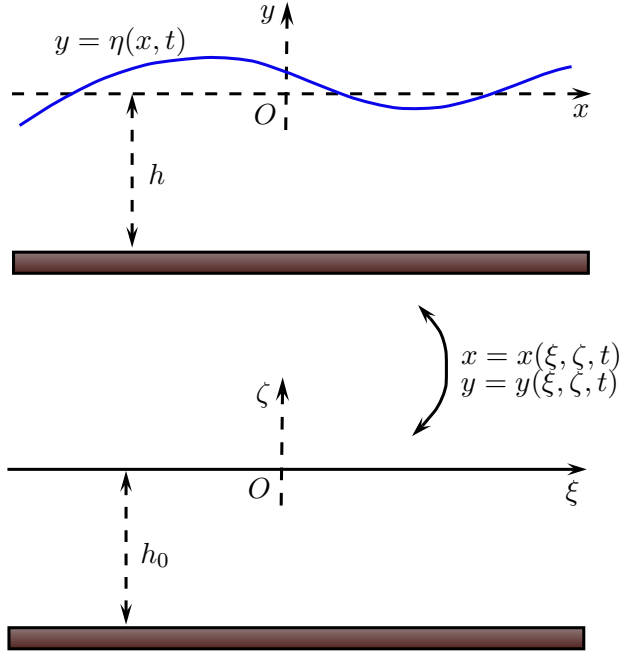
\begin{figure}
\centering
\scalebox{1} 
{
\begin{pspicture}(0,-4.3)(8.759063,4.3334374)
\definecolor{color191g}{rgb}{0.5490196078431373,0.4627450980392157,0.4588235294117647}
\definecolor{color191f}{rgb}{0.3176470588235294,0.13333333333333333,0.11764705882352941}
\definecolor{color194}{rgb}{0.043137254901960784,0.0196078431372549,0.9137254901960784}
\psframe[linewidth=0.04,dimen=outer,fillstyle=gradient,gradlines=2000,gradbegin=color191g,gradend=color191f,gradmidpoint=1.0](7.66,1.0)(0.06,0.68)
\psline[linewidth=0.04cm,linestyle=dashed,dash=0.16cm 0.16cm,arrowsize=0.043cm 3.0,arrowlength=1.6,arrowinset=0.4]{->}(0.06,3.08)(7.68,3.08)
\psline[linewidth=0.04cm,linestyle=dashed,dash=0.16cm 0.16cm,arrowsize=0.055cm 3.0,arrowlength=1.8,arrowinset=0.4]{<-}(3.68,4.3)(3.66,2.6)
\pscustom[linewidth=0.04,linecolor=color194]
{
\newpath
\moveto(0.16,2.66)
\lineto(0.5,2.88)
\curveto(0.67,2.99)(1.025,3.18)(1.21,3.26)
\curveto(1.395,3.34)(1.735,3.45)(1.89,3.48)
\curveto(2.045,3.51)(2.4,3.55)(2.6,3.56)
\curveto(2.8,3.57)(3.14,3.54)(3.28,3.5)
\curveto(3.42,3.46)(3.7,3.36)(3.84,3.3)
\curveto(3.98,3.24)(4.295,3.115)(4.47,3.05)
\curveto(4.645,2.985)(4.98,2.905)(5.14,2.89)
\curveto(5.3,2.875)(5.605,2.885)(5.75,2.91)
\curveto(5.895,2.935)(6.205,3.02)(6.37,3.08)
\curveto(6.535,3.14)(6.825,3.25)(6.95,3.3)
\curveto(7.075,3.35)(7.255,3.41)(7.42,3.44)
}
\psframe[linewidth=0.04,dimen=outer,fillstyle=gradient,gradlines=2000,gradbegin=color191g,gradend=color191f,gradmidpoint=1.0](7.6,-3.98)(0.0,-4.3)
\psline[linewidth=0.04cm,arrowsize=0.043cm 3.0,arrowlength=1.6,arrowinset=0.4]{->}(0.0,-1.9)(7.62,-1.9)
\psline[linewidth=0.04cm,linestyle=dashed,dash=0.16cm 0.16cm,arrowsize=0.055cm 3.0,arrowlength=1.8,arrowinset=0.4]{<-}(3.62,-0.68)(3.6,-2.38)
\usefont{T1}{ppl}{m}{n}
\rput(7.5445313,2.85){$x$}
\usefont{T1}{ppl}{m}{n}
\rput(3.4245312,4.13){$y$}
\usefont{T1}{ppl}{m}{n}
\rput(3.3945312,2.81){$O$}
\usefont{T1}{ppl}{m}{n}
\rput(7.434531,-2.21){$\xi$}
\usefont{T1}{ppl}{m}{n}
\rput(3.3645313,-0.91){$\zeta$}
\pscustom[linewidth=0.04]
{
\newpath
\moveto(5.46,-1.32)
\lineto(5.61,-1.14)
\curveto(5.685,-1.05)(5.765,-0.86)(5.77,-0.76)
\curveto(5.775,-0.66)(5.775,-0.465)(5.77,-0.37)
\curveto(5.765,-0.275)(5.72,-0.12)(5.68,-0.06)
\curveto(5.64,0.0)(5.59,0.07)(5.56,0.1)
}
\psline[linewidth=0.04cm,arrowsize=0.05291667cm 2.6,arrowlength=1.7,arrowinset=0.4]{->}(5.56,0.1)(5.4,0.28)
\psline[linewidth=0.04cm,arrowsize=0.055cm 2.7,arrowlength=1.7,arrowinset=0.4]{->}(5.46,-1.32)(5.32,-1.46)
\usefont{T1}{ppl}{m}{n}
\rput(7.0245314,-0.41){$x = x(\xi,\zeta,t)$}
\usefont{T1}{ppl}{m}{n}
\rput(7.0045314,-0.77){$y = y(\xi,\zeta,t)$}
\psline[linewidth=0.04cm,linestyle=dashed,dash=0.16cm 0.16cm,arrowsize=0.055cm 2.5,arrowlength=1.8,arrowinset=0.4]{<->}(1.64,3.08)(1.64,1.04)
\usefont{T1}{ppl}{m}{n}
\rput(1.9645313,2.01){$h$}
\psline[linewidth=0.04cm,linestyle=dashed,dash=0.16cm 0.16cm,arrowsize=0.055cm 2.5,arrowlength=1.8,arrowinset=0.4]{<->}(1.56,-1.92)(1.56,-3.96)
\usefont{T1}{ppl}{m}{n}
\rput(1.9445312,-3.03){$h_0$}
\usefont{T1}{ppl}{m}{n}
\rput(3.3145313,-2.15){$O$}
\usefont{T1}{ppl}{m}{n}
\rput(1.4945313,3.77){$y=\eta(x,t)$}
\end{pspicture} 
}
\caption{\small\em Conformal map of the physical domain into a uniform strip.}
\label{fig:cmap}
\end{figure}

The main idea behind handling the unknown free surface computationally is to reformulate the system using a time-dependent conformal map from the physical domain on a uniform strip. This transformation is schematically depicted in Figure~\ref{fig:cmap}. The new horizontal and vertical coordinates are denoted by $\xi$ and $\zeta$ respectively. The complete derivation of governing equations in the conformal domain can be found, \eg in \cite[Appendix A]{Mitsotakis2014a}. Here we just provide the final set of equations:
\begin{eqnarray}\label{eq:gam}
  \gamma_t &=& \gamma_\xi\,\H\Bigl[\frac{\psi_\xi}{J}\Bigr] - \chi_\xi\,\frac{\psi_\xi}{J}, \\
  \phis_t &=& \frac{1}{2}\frac{\psi_\xi^2 - \phis_\xi^2}{J} - g\gamma + \phis_\xi\,\H\Bigl[\frac{\psi_\xi}{J}\Bigr]. \label{eq:phis}
\end{eqnarray}
where $\gamma(\xi,t) := \eta(\chi(\xi,t),t)$, $\chi(\xi,t) := x(\xi,0,t)$ and the Jacobian $J := \chi_\xi^2 + \gamma_\xi^2$. Two evolution equations above have to be completed by two additional relations in order to close the system
\begin{equation*}
  \psi_\xi = \T[\phis_\xi], \qquad \chi_\xi = 1 - \H[\gamma_\xi],
\end{equation*}
where $\psi$ is the stream function (or equivalently the imaginary part of the complex potential). The pseudo-differential operators $\H$ and $\T$ are defined as
\begin{equation*}
  \hat{\H} = \ui\coth(kh_0), \qquad
  \H[f](x) = \frac{1}{2h_0}\int_{\R}f(y)\coth\Bigl(\frac{\pi}{2h_0}(y-x)\Bigr)\,\ud y,
\end{equation*}
\begin{equation*}
  \hat{\T} = \ui\tanh(kh_0), \qquad \T[f](x) = \frac{1}{2h_0}\int_{\R}f(y)\cosech\Bigl(\frac{\pi}{2h_0}(y-x)\Bigr)\, \ud y.
\end{equation*}

In order to discretize the free surface Euler system \eqref{eq:gam}, \eqref{eq:phis} one can use a Fourier-type pseudospectral method, where all derivatives along with nonlocal pseudo-differential operators are computed spectrally \cite{Boyd2000}. Nonlinear products are computed in real space and antialiased using the $3/2$ rule. The overall implementation is very efficient thanks to the FFT algorithm \cite{FFTW98}. As a result, we have an efficient solver for the full Euler equations with the free surface. The results of these simulations can provide the reference solution to validate and to assess the accuracy of various approximate models, see \eg \cite{Mitsotakis2014a}. The method described above was generalized recently to arbitrary bottoms in \cite{Viotti2014}. However, the derivation becomes significantly trickier and we do not enter into these details in this document.

\subsection{Computation of steady fully nonlinear solutions}

\begin{figure}
  \centering
  \subfigure[Horizontal velocity]{
  \includegraphics[width=0.48\textwidth]{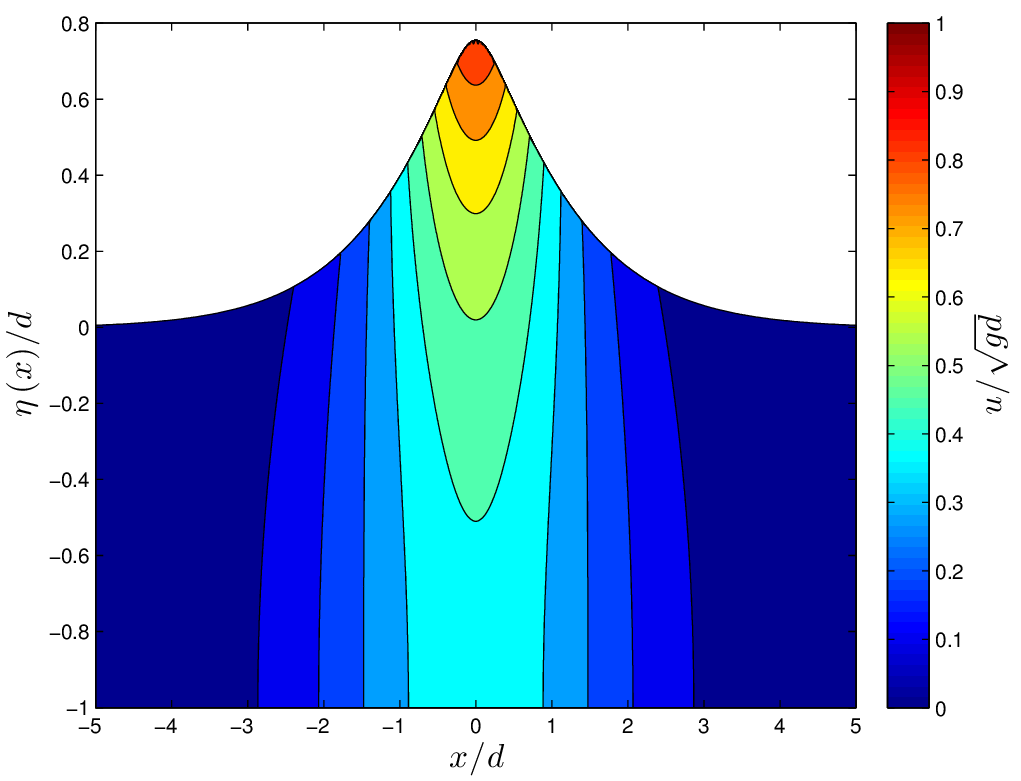}}
  \subfigure[Vertical velocity]{
  \includegraphics[width=0.48\textwidth]{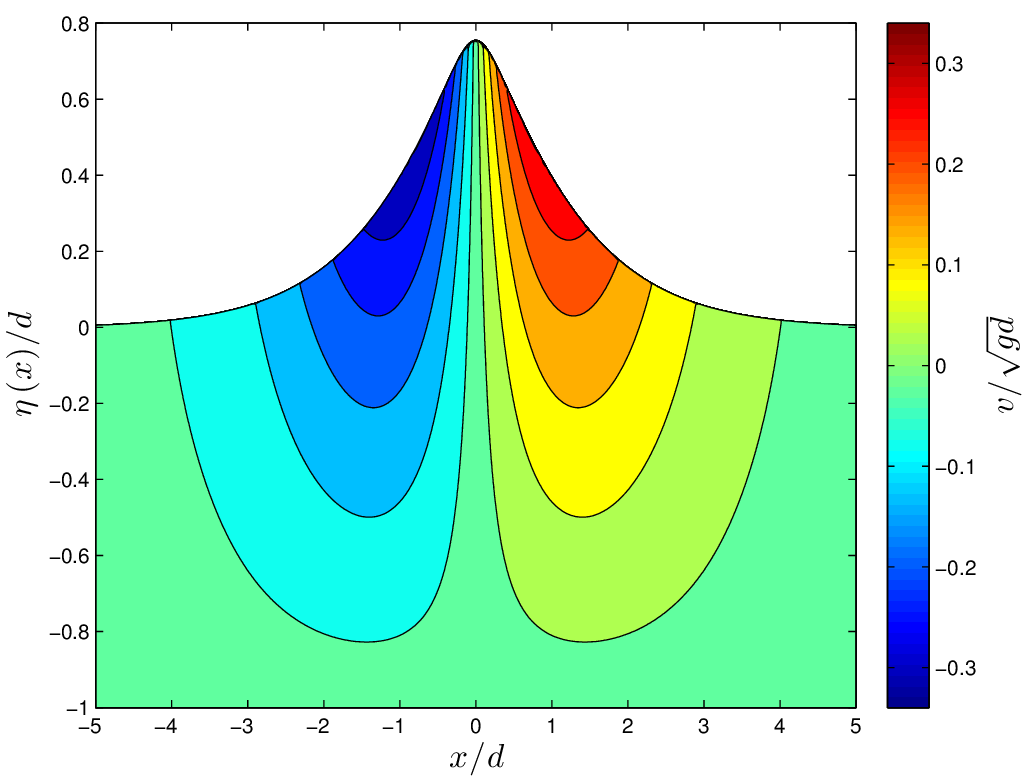}}
  \caption{\small\em Iso-horizontal (left) and iso-vertical (right) velocities under a large solitary wave. Lines correspond to the iso-values computed in the `fixed' Frame of reference where the the fluid is at rest in the far field $x\to\pm\infty$. Taken from \cite{Dutykh2013b}.}
  \label{fig:sw}
\end{figure}

The same conformal mapping technique can be used to compute efficiently steady travelling wave solutions to the full Euler equations with the free surface. Namely, we use the formulation found by K.~\textsc{Babenko} (1987) \cite{Babenko1987}, which has an advantage to have the same degree of nonlinearity (\ie quadratic for gravity waves) as the original Euler equations. See \cite{Clamond2012b, Dutykh2013b} for more details on the derivation of Babenko's equation for gravity solitary waves. The resulting Babenko's equation is solved numerically using the iterative Petviashvili scheme \cite{Petviashvili1976}. This method can be implemented in a very elegant way within fifty lines of \texttt{Matlab} code (for the computational core). The code  is freely available in Open Source at the \texttt{Matlab File Exchange} server \cite{Clamond2012}. Another advantage of this formulation is that it allows to compute physical fields in the fluid bulk as well (not only at the free surface). For the sake of illustration we show the velocity field on Figure~\ref{fig:sw}. This formulation has been recently extended for the computation of capillary--gravity waves as well \cite{Clamond2015a}. As above, these numerically exact solutions can be used to validate various approximate long wave models which have to reproduce the corresponding solitary wave solutions of the base model with high physical fidelity.

\section{Extended fully nonlinear shallow water equations}

For two-dimensional surface water waves propagating in shallow water of constant depth, one can approximate the velocity field by \cite{Dutykh2014c}
\begin{equation*}
  u(x,y,t)\ \approx\ \bar{u}(x,t), \qquad v(x,y,t)\ \approx\ -\,(y+d)\,\bar{u}_x
\end{equation*}
where $\,d\,$ is the mean water depth and $\,\bar{u}\,$ is the horizontal velocity averaged over the water column --- \ie $\,\bar{u}\equiv h^{-1}\int_{-d}^\eta u\,\ud y$ --- $\,y=\eta\,$ and $\,y=0\,$ being the equations of the free surface and of the still water level, respectively. The horizontal velocity $\,u\,$ is thus uniform along the water column and the vertical velocity $\,v\,$ is chosen so that the fluid incompressibility if fulfilled. \textsc{Serre} (1953) \cite{Serre1953} derived the following approximate system of equations
\begin{align}
h_t\ +\ \partial_x\!\left[\,h\,\bar{u}\,\right]\, &=\ 0, \label{eq:massse} \\
\partial_t\!\left[\,h\,\bar{u}\,\right]\, +\ \partial_x\!\left[\,h\,\bar{u}^2\,
+\,\half\,g\,h^2\,+\,\third\,h^2\,\gamma\,\right]\, &=\ 0, \label{eq:qdmfluxse} 
\end{align}
where 
\begin{equation}\label{eq:defaccver}
\gamma\ =\ h\,(\bar{u}_x^{\,2}-\bar{u}_{xt}-\bar{u}\/\bar{u}_{xx})\ =\ 2\,h\,\bar{u}_x^{\,2}\ -\ h\,\partial_x\!\left[\,\bar{u}_t\,+\,\bar{u}\,\bar{u}_x\,\right],
\end{equation}
is the vertical acceleration of the fluid at the free surface \cite{Clamond2009}. Physically, equations \eqref{eq:massse} and \eqref{eq:qdmfluxse} describe, respectively, the mass and momentum flux conservations.

It is however possible to obtain a more general system by intoducing a free parameter into the model using Bona--Smith--Nwogu's trick \cite{BS, Nwogu1993}. Here we provide only the final result, while the computations can be found in \cite{Dutykh2014c}:
\begin{align}
  h_t\ +\ \partial_x\!\left[\,h\,\bar{u}\,\right]\, &=\ 0,  \label{eq:masssem} \\
  \partial_t\!\left[\,h\,\bar{u}\,\right]\, +\ \partial_x\!\left[\,h\,\bar{u}^2\, +\,\half\,g\,h^2\,+\,\third\,h^2\,\gamma\,\right]\, &=\ 0, \label{eq:qdmfluxsem} \\
  2\,h\,\bar{u}_x^{\,2}\ +\ (1-\alpha)\,g\,h\,h_{xx}\ -\ \alpha\,\,h\,\partial_x\!\left[\,\bar{u}_t\,+\,\bar{u}\,\bar{u}_x\,\right]\,&=\ \gamma,
\label{eq:veraccmod}
\end{align}
where $\alpha$ is a constant at our disposal. We note however that we were not able to find the energy conservation equation for the last system, even if it is Galilean invariant \cite{Dutykh2014c}.

\subsection{Linear approximation}

For infinitesimal waves, $\eta$ and $\bar{u}$ being both small, it is reasonable to linearise the equations around $\eta = 0$ and $\bar{u} = 0$. We obtain thus the linear system of equations
\begin{align}
  \eta_t\ +\ d\,\bar{u}_x\ &=\ 0,  \label{eqmasssemlin} \\
  \bar{u}_t\ +\ g\,\eta_x\ +\ \third\,d\,\gamma_x\ &=\ 0, \label{eqqdmfluxsemlin} \\
  (1-\alpha)\,g\,d\,\eta_{xx}\ -\ \alpha\,\,d\,\bar{u}_{xt}\ &=\ \gamma. \label{eqveraccmod}
\end{align}
Seeking for traveling waves of the form $\eta=a\cos \bigl(k(x-ct)\bigr)$, we obtain the (linear) dispersion relation
\begin{equation}\label{eq:disrelserlin}
  \frac{c^2}{g\/d}\ =\ \frac{3\,+\,(\alpha-1)\/(k\/d)^2}{3\,+\,\alpha\/(k\/d)^2}\ =\ 1\, -\, \third\/(k\/d)^2\, +\,\textstyle{1\over9}\/\alpha\/(k\/d)^4\,-\, \textstyle{1\over27}\/\alpha^2\/(k\/d)^6\,+\,\cdots.
\end{equation}
We note that this relation is well-posed (\ie $c^2>0$ for all $k$) only if $\alpha \geqslant 1$. In order to find a suitable choice for $\alpha$, the relation \eqref{eq:disrelserlin} can be compared with the dispersion relation of linear waves on finite depth
\begin{equation*}
  c^2\,/\,g\,d\ =\ \mathrm{thc}(k\/d)\ =\ 1\, -\, \third\/(k\/d)^2\, +\,\textstyle{2\over15}\/(k\/d)^4\,-\,\textstyle{17\over315}\/(k\/d)^6\, + \,\textstyle{62\over2835}\/(k\/d)^8\, + \,\cdots,
\end{equation*}
where $\mathrm{thc}(x)\equiv\tanh(x)/x$ if $x \neq 0$ and $\mathrm{thc}(0) \equiv 1$. Comparing the Taylor expansions, it is clear that \eqref{eq:disrelserlin} matches the exact one only up to the second-order in general, except when $\alpha=6/5$ in which case it matches up to the fourth-order. Therefore, $\alpha_{\mathrm{opt}} = 6/5$ is a suitable choice having the advantage of being independent of the wave characteristics. This method of choosing the optimal $\alpha$ has been used by many authors starting from the pioneering works \cite{BS, Nwogu1993, Beji1996}.

\subsection{Validation on solitary waves}

\begin{figure}
  \centering
  \includegraphics[width=0.99\textwidth]{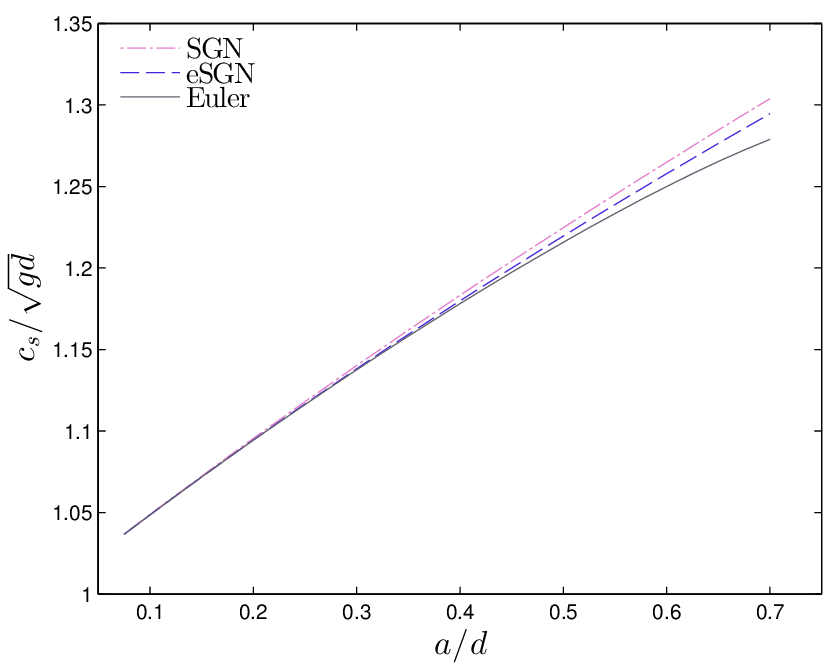}
  \caption{\small\em Speed--amplitude relations for solitary waves in \acs{SGN}, \acs{eSGN} and the full Euler equations ($\alpha = 6/5$).}
  \label{fig:spa}
\end{figure}

We will compare the solitary wave solutions to the three following models:
\begin{itemize}
  \item \acs{SGN} equations
  \item \acs{eSGN} equations (with optimal $\alpha$)
  \item the full Euler equations (the reference solution)
\end{itemize}
The solitary waves to the classical \acs{SGN} equations are known analytically. The solitary wave solutions for the full Euler equations are computed using the method of conformal variables \cite{Clamond2012b, Dutykh2013b}. The \textsc{Matlab} script used to generate the solitary waves can be downloaded at \cite{Clamond2012}. Unfortunately, we did not succeed in finding analytical solutions to the \acs{eSGN} equations for a general $\alpha$. Consequently, we had to employ the numerical methods. The propagation speeds predicted by various models are reported in Table~\ref{tab:cs}. One can already see that the \acs{eSGN} predictions are always closer to the full Euler equations. We computed the speed--amplitude relation for the whole range of amplitudes (see Figure~\ref{fig:spa}). One can see that the \acs{eSGN} model approximates better the reference solution again. It is somehow surprising, since the model was tuned on linear solution and this ``tuning'' turns out to improve the description of nonlinear solutions as well. Currently, we are looking for a model which does possess the full set of conservation laws together with improved dispersion relation properties.

\begin{table}
  \centering
  \begin{tabular}{||>{\columncolor[gray]{0.85}}l||>{\columncolor[gray]{0.85}}c||>{\columncolor[gray]{0.85}}c||>{\columncolor[gray]{0.85}}c||}
  \hline\hline
  \textit{Model/Amplitude} & $a/d = 0.1$ & $a/d = 0.45$ & $a/d = 0.70$ \\
  \hline\hline
  \acs{SGN} speed,  $c_s/\sqrt{gd}$ & \textbf{1.048}80  & \textbf{1.2}041 & \textbf{1.}3038 \\ 
  \acs{eSGN} speed, $c_s/\sqrt{gd}$ & \textbf{1.0485}6  & \textbf{1.19}99 & \textbf{1.2}946 \\ 
  Full Euler speed, $c_s/\sqrt{gd}$ & \textbf{1.048548} & \textbf{1.1973} & \textbf{1.2788} \\ 
  \hline\hline
  \end{tabular}
  \bigskip
  \caption{\small\em Comparison of the solitary wave speeds for several fixed values of the wave amplitude. The parameter $\alpha = 6/5$.}
  \label{tab:cs}
\end{table}

\appendix
\section{Shallow water waves on rough bottoms}

\begin{figure}
  \centering
  \includegraphics[width=0.95\textwidth]{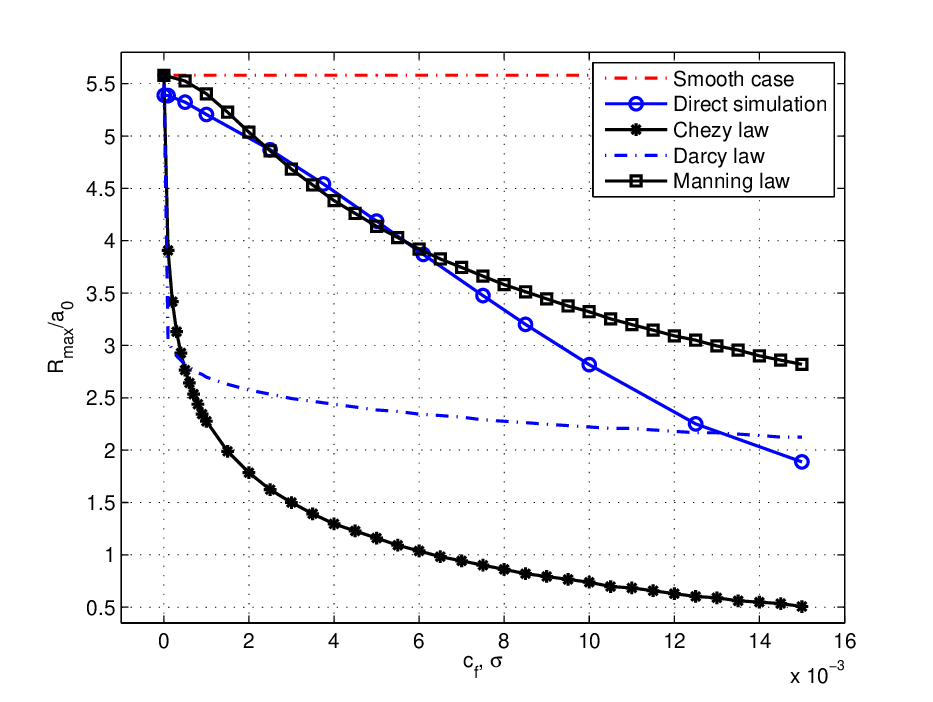}
  \caption{\small\em Comparison of the runup reduction effect for various ad-hoc friction terms and the random bottom perturbation model. The horizontal axis represents the friction coefficient $c_f$ for deterministic computations and $\sigma$ for the random roughness model (blue solid line).}
  \label{fig:runup}
\end{figure}

In this final section we would like to mention a result, which appeared in a Physics Journal (PRL) \cite{Dutykh2011c} and it remains essentially unknown in the hyperbolic community. The authors considered the wave run-up problem on random rough bottoms. Our stochastic computations were compared to several simulations using classical friction terms routinely used to model the bottom rugosity. Surprisingly, a very good qualitative agreement was obtained with the Manning--Strickler law, while the Ch\'ezy and Darcy--Weisbach laws provide too strong momentum damping. See Figure~\ref{fig:runup} for the numerical results and \cite{Dutykh2011c} for more details. As a conclusion, we recommend to use the Manning--Strickler law to take into account the wave propagation on rough bottoms.

\backmatter
\bibliography{biblio}
\bibliographystyle{smfplain}

\end{document}